\documentclass [12pt,english] {article}
\usepackage[utf8]{inputenc}
\usepackage{amsthm}
\usepackage{amsmath}
\usepackage{amsfonts}
\usepackage{amssymb}
\usepackage[all]{xy}
\usepackage{color}   
\usepackage{tikz} \usetikzlibrary{arrows}
\usetikzlibrary{shapes.geometric}
\usepackage{listings}
\usepackage{setspace}
\usepackage[font=small,labelfont=bf]{caption}
\singlespace

\newtheorem{ejem}{Example}
\newtheorem{defi}[ejem]{Definition}
\newtheorem{teo}[ejem]{Theorem}
\newtheorem{prop}[ejem]{Proposition}
\newtheorem{lema}[ejem]{Lemma}
\newtheorem{remark}[ejem]{Remark}
\newtheorem{coro}[ejem]{Corollary}
\newtheorem{conje}[ejem]{Question}
\numberwithin{ejem}{section}

\usepackage[textwidth=3cm, textsize=small, colorinlistoftodos]
{todonotes}

\title{Generalized Token Graphs}
\author{C. Amairani Herrera-Ramirez \footnote{Universidad Veracruzana, Xalapa, M\'exico, citlaliamairaniherreraramirez@gmail.com}, Teresa I. Hoekstra-Mendoza \footnote{ Centro de Investigaci\'on en Matem\'aticas, Guanajuato, M\'exico, maria.idskjen@cimat.mx, ORCID:0000-0002-8012-6054, (Corresponding author)} }

\begin{document}
\maketitle
\begin{abstract}
    In this paper we give a new generalization of token graphs. Given two integers $1\leq m \leq k$ and a graph $G$ we define the generalized token graph of the graph $G$, to be the graph  $F_k^m(G)$ whose vertices correspond to configurations of $k$ indistinguishable tokens placed at distinct vertices of $G$, where  two configurations are adjacent whenever one configuration can be reached from the other by moving $m$ tokens along $m$ edges of $G$. When $m=1$, the usual token graph $F_k(G)$ is recovered. 
    We give sufficient and necessary conditions on the graph $G$ for $F_2^2(G)$ to be connected and
    we give sufficient and necessary conditions on the graph $G$ for $F_2^2(G)$ to be bipartite. 
    We also analyse some properties of generalized token graphs, such as  clique number, chromatic number, independence number and domination number.  Finally, we conclude with an analysis of the automorphism group of the generalized token graph.
\end{abstract}

MSC: 05C76, 05C40, 05C15 

Keywords: Generalized Token Graphs, Connectivity, Bipartiteness, Automorphism group

\section*{Acknowledgments}
The author Teresa I. Hoekstra-Mendoza was supported by Sistema Nacional de Investigadores (SNI) number 829061 which is gratefully acknowledged.

\section{Introduction}
Many problems in mathematics are modelled by moving objects on the vertices of a graph according to certain prescribed rules. In ordinary token graphs, $k$ indistinguishable objects ocupy $k$ different vertices and they are allowed to move one at a time. Thus a natural generalization is to allow a fixed number of tokens to move at once. This motivates the definition of generalized token graph.

\begin{defi}
Let $G$ be a graph and $m \leq k$ two positive integers. We define the generalized token graph $F_k^m(G)$ of $G$ as $V(F_k^m(G))=\{B\subset V(G): |A|=k\}$ and, if $A=\{a_1, \dots, a_k\}$, $B=\{b_1, \dots, b_k\}$ then $(A,B)\in E(F_k^m(G))$ if there are $m$ indices $i_1, \dots, i_m$ such that $(a_{i_j}, b_{i_j})\in E(G)$ for $ 1\leq j \leq m$ and $a_l=b_l$ for every $l \neq i_j $.
\end{defi}

Generalized token graphs, for example, can model the movement of trains. Consider a graph where the vertices correspond to train stations and two of them are adjacent if there exists a train track between them which does not pass through any other station. If only one train can move at a time, the transport system would be too slow, so several trains have to move at the same time without collisions.

\begin{ejem}
    Consider the cycle on four vertices $C_4=K_{2,2}$. Then the graph $F_2^2(C_4)$ is  isomorphic to the disjoint union $K_4\sqcup K_2.$
\end{ejem}

In \cite{R}, the following generalizations of token graphs were mentioned as an open line of investigation for future work:
\begin{defi}\cite{R}
For $r \in \{1,\dots, k\}$, let $F_{k,r}(G)$ be the graph with
vertex set $\binom{V(G)}{k},$ where two vertices $A$ and $B$ in $F_{k,r}(G)$ are adjacent whenever $|A\Delta B| = 2r$
and there is a perfect matching between $A \setminus B$ and $B \setminus A$ in $G$.
Let $F'_{k,r}(G)$ be the variant
where instead we require that every edge is present between $A \setminus B$ and $B \setminus A.$
\end{defi}
These two generalizations are also generalizations of the Kneser graph, $KG_{n,k}$, whose vertices are the $k$-subsets of an $n$-set, where two vertices $A$ and $B$ are adjacent
whenever $A \cap  B = \emptyset$. However, both generalizations $F_{k,m}(G)$ and $F'_{k,m}(G)$ are different from our generalization $F_k^m(G)$, since in both cases, for two vertices $A,B$ to be adjacent in $F_{k,m}(G)$ or $F'_{k,m}(G)$, it is necessary that, the symmetric difference $A\Delta B$ is exactly $2m$ while in $F_k^m(G) $ we allow $|A\Delta B|$ to be less or equal than $2m$.

\begin{prop}Let $x,y\in V(F_2^2(G))$. Assume that $v$ and $w$ have $k$ common neighbours. If $(v,w)\notin E(G)$ then $d(vw)=d(v)d(w)-\frac{k}{2}(k+1)$ and if $(v,w)\in E(G)$ then $d(vw)=d(v)d(w)-\frac{k}{2}(k+1)-1$.
    
\end{prop}
\begin{proof}
   Assume first that $(v,w)\notin E(G)$. Notice that, for a neighbour $xy$ of $vw$ there are four options:
   \begin{itemize}
       \item $x,y \in N(v)\cap N(w)$, in which case there are $\frac{k(k-1)}{2}$ options for $xy$
       \item $x,y \notin N(v)\cap N(w)$, in which case there are $(d(w)-k)(d(v)-k)$ options for $xy$
       \item $x\in N(v)\cap N(w), y \in N(w)\setminus(N(v)\cap N(w))$ in which case there are $k(d(w)-k)$ options for $xy$, and 
       \item $x\in N(v)\cap N(w), y \in N(v)\setminus(N(v)\cap N(w))$ in which case there are $k(d(v)-k)$ options for $xy$.
   \end{itemize}
   This means that the total number of options for $xy$ is $\frac{k(k-1)}{2}+(d(w)-k)(d(v)-k)+k(d(w)-k)+k(d(v)-k)=d(v)d(w)-k^2+\frac{1}{2}k(k-1)=d(v)d(w)-\frac{k}{2}(k+1).$

   \begin{center}
		\begin{tikzpicture}[scale=2]
			\draw (-2.69,1.95) node[draw,circle, scale=.5] (0) { v};
			\draw (-0.61,1.97) node[draw,circle, scale=.5] (1) { w};
			\draw (-1.6,1.62) node[draw,circle, scale=.5] (2) { };
			\draw (-1.61,1.33) node[draw,circle, scale=.5] (3) { };
			\draw (-1.61,0.85) node[draw,circle, scale=.5] (4) { };
			\draw (-1.62,1.25) node[draw, ellipse, minimum width=20pt, minimum height=70pt] (5) { };
			\draw (-0.63,1.35) node[draw,circle, scale=.5] (6) { };
			\draw (-0.37,1.35) node[draw,circle, scale=.5] (7) {};
			\draw (0.25,1.35) node[draw,circle, scale=.5] (8) {};
			\draw (-0.2,1.35) node[draw, ellipse, minimum width=75pt, minimum height=20pt] (9) {};
            \draw (-0.2,1) node (w) {$d(w)-k$};
            \draw (-3.2,1) node (v) {$d(v)-k$};
			\draw (-3.53,1.35) node[draw,circle, scale=.5] (10) { };
			\draw (-3.21,1.35) node[draw,circle, scale=.5] (11) {};
			\draw (-2.69,1.35) node[draw,circle, scale=.5] (12) { };
            \draw (-1.3,0.8) node (x) {$k$};
			\node[draw, ellipse, minimum width=75pt, minimum height=20pt] (13) at (-3.1,1.35)  { };
			\draw  (1) edge (2);
			\draw  (0) edge (2);
			\draw  (1) edge (3);
			\draw  (1) edge (4);
			\draw  (0) edge (3);
			\draw  (0) edge (4);
			\draw[dotted]  (3) edge (4);
			\draw[dashed]  (0) edge (1);
			\draw[dotted]  (7) edge (8);
			\draw  (1) edge (6);
			\draw  (1) edge (7);
			\draw  (1) edge (8);
			\draw  (0) edge (12);
			\draw  (0) edge (11);
			\draw  (0) edge (10);
			\draw[dotted]  (11) edge (12);
		\end{tikzpicture}
	\end{center}

The case when  $(v,w)\in E(G)$ is analogous.
\end{proof}


\section{Main results}
One of the main differences between usual token graphs and generalized token graphs, is that when $G$ is connected, the usual token graph $F_k(G)$ is always connected. However, this is not the case for generalized token graphs.

\begin{teo}
    Let $G$ denote a connected bipartite graph. Then $F_2^2(G)$ is disconnected.
\end{teo}
\begin{proof}
    Let $V(G)=X\cup Y$ where both $X$ and $Y$ are independent sets. Consider the sets $Z=\{xy\in V(F_2^2(G)): x\in X, y\in Y\}$ and $W =\{vw \in V(F_2^2(G)): v,w \in X, \text{ or } v,w\in Y\} =  V(F_2^2(G))\setminus Z.$ Then there are no edges between $Z$ and $W$, thus $ F_2^2(G)$ is disconnected.
\end{proof}

\begin{teo}
    The graph $F_2^2(K_{1,n})$ consists of a complete subgraph $K_n$, together with $\binom{n}{2}$ isolated vertices.
\end{teo}
\begin{proof}
    Denote by $x_0$ the unique vertex of degree $n.$ Then, the subgraph generated by all configurations of the form $\{x_0y\}$ forms a complete subgraph on $n$ vertices. Now consider a configuration $\{yz\}$ which does not contain the vertex $x_0.$ Since in $G$, we have $N(y)=N(z)=\{x_0\}$, $\{yz\}$ corresponds to an isolated vertex and there are $\binom{n}{2}$ such configurations.
\end{proof}
\begin{coro}
    The independence number and domination number of $F_2^2(K_{1,n})$ are $$\gamma(F_2^2(K_{1,n}))= \binom{n}{2}+1 = \alpha (F_2^2(K_{1,n}))$$
\end{coro}
\begin{teo}\label{kmn}
Let $K_{m,n}$ denote the complete bipartite graph with $m,n>1$. Then $F_2^2(K_{m,n})$ consists of the disjoint union of a complete subgraph $K_{mn}$ with a complete bipartite subgraph $K_{\binom{n}{2},\binom{m}{2}}.$
\end{teo}

\begin{proof}
    Let $M$ and $N$ be the independent subsets of $V(K_{m,n})$ or cardinalities $m$ and $n$ respectively. Then any vertex of the form $xy$ with $x \in M$ and $y \in N$ is adjacent to any vertex $zw$ with $z\in M$ and $w \in N$ since $(x,w), (y,z)\in E(K_{m,n}).$ This means that the subgraph induced by the vertices $\{xy: x\in M, y\in N\}$ is isomorphic to $K_{mn}.$ Notice that $\{ab: a,b\in N\}$ and $\{uv: u,v\in M\}$ are independent subsets in $F_2^2(K_{m,n})$ since both $M$ and $N$ are independent sets. Now, between $ab$ and $uv$ with $a,b\in N, uv\in M$ there exists every edge since $(a,u),(a,v),(b,u),(b,v)\in E(K_{m,n})$ thus we obtain a subgraph isomorphic to $K_{\binom{n}{2},\binom{m}{2}}.$ Finally, notice that between these two subgraphs there are no edges.
\end{proof} 

\begin{coro}
   For $n,m\geq 2,$ the chromatic number, clique number, domination number and independence number for $F_2^2(K_{m,n})$  are the following $$\chi (F_2^2(K_{m,n}))=mn= \omega(F_2^2(K_{mn})), \gamma(F_2^2(K_{m,n}))=3, \text{ and }$$  $$\alpha(F_2^2(K_{m,n}))=max\{\binom{n}{2},\binom{m}{2} \}+1.$$
\end{coro}

For the general case where we allow all tokens to move at the same time on a complete bipartite graph, we have the following result.

\begin{teo}
Let $G=K_{m,n}$ be the complete bipartite graph for $m,n>1$ and $k \leq min\{m,n\}$. Then $$F_k^k(K_{m,n})=K_{\binom{n}{0}\binom{m}{k},\binom{n}{k}\binom{m}{0}} \sqcup K_{\binom{n}{k-1}\binom{m}{1},\binom{n}{1}\binom{m}{k-1}} \sqcup \dots \sqcup H,$$ where $$H=\left\{\begin{array}{cc}
    K_{\binom{n}{\frac{k}{2}}\binom{m}{\frac{k}{2}}} & \mbox{if }k\equiv 0\mod2\\
    K_{\binom{n}{k-\lfloor\frac{k}{2}\rfloor}\binom{m}{\lfloor\frac{k}{2}\rfloor},\binom{n}{\lfloor\frac{k}{2}\rfloor}\binom{m}{k-\lfloor\frac{k}{2}\rfloor}} & \mbox{if }k\equiv 1\mod2
\end{array}\right.$$
\end{teo}
\begin{proof}
Let $V(K_{m,n})=M\cup N$ with $|M|=m$ and $|N|=n$.
Given integers $p>q \geq 0$ such that $k=p+q$,  let $A_{p,q}$ denote the set of token configurations which have $p$ tokens on vertices of $M$ and $q$ tokens on vertices of $N$ and notice that  $|A|=\binom{m}{p}\binom{n}{q}$. Let $B_{q,p}$ denote the set of token configurations which have $ p$ tokens on vertices in $N$ and $q$ tokens on vertices of $M$ and notice that $|B|=\binom{m}{q}\binom{n}{p}$. Any vertex of $B$ is adjacent to every vertex in $A$ and any vertex of $A$ is adjacent to every vertex of $B$. If $p\neq q$ then both $A_{p,q}$ and $B_{q,p}$ are independent sets.
If $k=2l$, we have $\{p,q \in \mathbb{Z}: p+q=k\}=\{(1,k-1),(2,k-2),(3,k-3),\dots,(l,l)\}$, and if $k=2l+1$ we have $\{p,q \in \mathbb{Z}: p+q=k\}=\{(1,k-1),(2,k-2),(3,k-3),\dots,(l,l+1)\}$.
Thus, for any integers $p,q\in \mathbb{Z}^+\cup \{0\}$ such that $p+q=k$ and $p\neq q$ we have a connected component in $F_k^k(K_{m,n})$ which is isomorphic to $$K_{\binom{m}{p}\binom{n}{q},\binom{m}{q}\binom{n}{p}}.$$
If $p=  \frac{k}{2} =q$  we have $A_{p,q} = B_{q,p}$, thus we obtain a connected component isomorphic to $K_t$ where $t=\binom{n}{\frac{k}{2}}\binom{m}{\frac{k}{2}}.$



\end{proof}

We shall now analyse the generalized token graph of the cycle graph $C_n$, which depends strongly of the parity of $n$, and of the congruency of $n$ modulo 4.

\begin{teo}\label{prod}
    The graph $F_2^2(C_n)$ is isomorphic to the graph $C_n \square P_{\frac{n-1}{2}}$ for any odd number $n.$
\end{teo}

\begin{proof}
Assume that $V(C_n)=\{1, \dots, n\}$.
    Consider the sets $V_i=\{(x,y)\in V(F_2^2(C_n)): d(x,y)=i\}$ and notice that $$V(F_2^2(C_n))=\bigcup_{i=1}^{\frac{n-1}{2}}V_i.$$ Every subgraph spanned by a set $V_i$ is a single cycle of length $n.$ The edges of $F_2^2(C_n)$ between different sets $V_i$ occur only between $V_i$ and $V_{i\pm 2}, $ to be specific, for every $(x,y)\in V_i$ we have the edges $(x-1,y+1), (x+1,y-1)\in E(F_2^2(C_n))$ and every other neighbour of $(x,y)$ belongs to $V_i.$ Hence, $F_2^2(C_n)\cong C_n\square P_{\frac{n-1}{2}}.$
\end{proof}

\begin{coro}
    For any odd number $n$, $\alpha(F_2^2(C_n))=(\frac{n-1}{2})^2.$
\end{coro}

In \cite{NPA}, the following bounds for the domination number of cylindrical grid graphs were proven.

\begin{teo} \cite{NPA}
  For $m,n\in \mathbb{Z},$ $$\frac{mn}{5}\leq \gamma(P_n\square C_m) \leq \frac{(m+2)(n+2)}{5}.$$  
\end{teo}
\begin{coro}
    The domination number of the generalized token graph is bounded by 
    $$\frac{m(2m+1)}{5}\leq F_2^2(C_{2m+1})\leq \frac{(m+2)(2m+3)}{5}.$$
\end{coro}

Moreover, in \cite{NPA}, the following exact values for the domination numbers of cylindrical grid graphs were obtained.

\begin{teo} \cite{NPA}
    The following domination numbers are 
    $$\gamma(C_3\square P_1) = 1, \gamma(C_5\square P_2)=3, \gamma(C_7\square P_3)=6, \gamma(C_9\square P_4) \text{ and } \gamma(C_{11}\square P_5)=14 $$
\end{teo}

\begin{coro}
    For $ m \leq 5,$ the exact values of $\gamma(F_2^2(C_{2m+1}))$ are $$\gamma (F_2^(C_3)) = 1, \gamma(F_2^2(C_5))=3, \gamma(F_2^2(C_7))=6, \gamma(F_2^2 (C_9)) \text{ and } \gamma(F_2^2(C_{11}))=14 .$$
\end{coro}

Another significant difference between usual token graphs and generalized token graphs is the following.
In usual token  graphs, if a graph $G$ is bipartite, then its token graph $F_k(G)$ is always bipartite, however this is not true for generalized token graphs.

\begin{prop}\label{comp}
    The graph   $F_2^2(C_n)$ contains a bipartite connected component and a non-bipartite connected component for any even number $n>4.$
\end{prop}
\begin{proof}
Assume that $V(C_n)=\{1, \dots, n\}$.
     The vertices $\{1,n\},\{2,n+1\}, \dots, \{n+1,2n\}$ form a cycle of length $n+1$ which is odd. This means that $F_2^2(C_n)$ is not bipartite. Notice that the set $\{(2x,2y), (2x+1,2y+1)\in V(F_2^2(C_n)): 1\leq x,y\leq \frac{n}{2} \}$ forms a connected component. Moreover, this connected component is bipartite since the sets $\{(2x,2y)\in V(F_2^2(C_n)): 1\leq x,y\leq \frac{n}{2} \}$ and $\{(2x+1,2y+1)\in V(F_2^2(C_n)): 1\leq x,y\leq \frac{n}{2} \}$ are  both independent. Hence, the odd cycle is contained in the connected component given by the vertices $\{(2x,2y+1)\in V(F_2^2(C_n)): 1\leq x,y\leq \frac{n}{2} \}$.
\end{proof}

\begin{figure}[h!]
\begin{center}
		\begin{tikzpicture}[scale=1,thick]
		\tikzstyle{every node}=[minimum width=0pt, inner sep=2pt, circle]
			\draw (-1.32,0.02) node[draw] (0) {};
			\draw (0.67,-1.14) node[draw] (1) {};
			\draw (0.67,1.2) node[draw] (2) {};
			\draw (-2,0.01) node[draw] (3) {};
			\draw (-1.0,-1.73) node[draw] (4) {};
			\draw (1.0,-1.73) node[draw] (5) {};
			\draw (2.0,-0.01) node[draw] (6) {};
			\draw (1,1.71) node[draw] (7) { };
			\draw (-1,1.72) node[draw] (8) {};
			\draw (4,0) node[draw] (9) {};
			\draw (5,-1.732050807568877) node[draw] (10) {};
			\draw (7,-1.7320508075688776) node[draw] (11) {};
			\draw (8,0) node[draw] (12) {};
			\draw (7,1.7) node[draw] (13) {};
			\draw (5,1.7320508075688772) node[draw] (14) {};
			\draw  (0) edge (1);
			\draw  (0) edge (2);
			\draw  (1) edge (2);
			\draw  (3) edge (4);
			\draw  (4) edge (5);
			\draw  (5) edge (6);
			\draw  (6) edge (7);
			\draw  (7) edge (8);
			\draw  (3) edge (8);
			\draw  (2) edge (7);
			\draw  (2) edge (4);
			\draw  (1) edge (5);
			\draw  (1) edge (8);
			\draw  (0) edge (3);
			\draw  (0) edge (6);
			\draw  (9) edge (10);
			\draw  (10) edge (11);
			\draw  (11) edge (12);
			\draw  (12) edge (13);
			\draw  (13) edge (14);
			\draw  (9) edge (14);
			\draw  (10) edge (13);
			\draw  (9) edge (12);
			\draw  (11) edge (14);
		\end{tikzpicture} 
        \caption{The graph $F_2^2(C_6)$.} \label{C6}
	\end{center}
\end{figure}
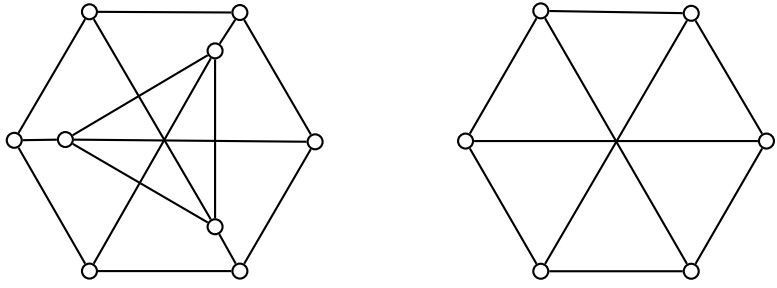

\begin{remark}\label{nonbi}
    The non-bipartite component of the graph $F_2^2(C_{2k})$ has the following description. Assume that $V(C_{2k})=\{1, \dots, n=2k\}$.
    For any odd number $i <  \frac{n}{2}$, define the sets $V_i=\{(x,y)\in V(F_2^2(C_n)): d(x,y)=i\}$ and let $k'$ denote the greatest odd integer such that $k'< \frac{n}{2}$.
Notice that $V_1\cup V_3 \cup \dots \cup V_{k'}$ is the vertex set of the non-bipartite connected component of $F_2^2(C_n)$.
 For $i<k'$, the subgraph spanned by a set $V_i$ is a cycle of length $n.$ If $n\equiv 0 (\text{mod }4)$ then the subgraph spanned by $V_{k'}$ is the graph obtained by a cycle of length $n$ by adding the edges $((x,y), (x+2k,y+2k))mod$ $n$. If $n\equiv 2 (\text{mod }4)$ then the subgraph spanned by $V_{k'}$ is a cycle of odd length $\frac{n}{2}.$  Given a vertex $(x,y)\in V_i$, its only neighbours outside of $V_i $ are $(x+1,y-1)$ and $(x-1,y+1)$.
\end{remark}

\begin{ejem}
    In Figure \ref{C6} we can see both connected components of the graph $F_2^2(C_6).$
\end{ejem}

\begin{teo}
    The chromatic numbers and clique numbers of the graph $F_2^2(C_n)$ are $ \chi(F_2^2(C_n))=3$ if $n \neq 4$ and $\omega (F_2^2(C_n))= 2$ if $n\neq 3,4,6.$ Moreover, $\omega (F_2^2(C_n))= 3$ for $n=3,6$ and $\omega (F_2^2(C_4))= 4=\chi((F_2^2(C_4)).$
\end{teo}

\begin{proof}
We start by analyzing the clique number.
Assume first that $F_2^2(C_n)$ contains a triangle, and notice that there are four possible cases, where $a,b,c,d,e,f$ are all distinct vertices of $C_n.$
\begin{center}
    \begin{tikzpicture}[scale=1.2]
        \node[circle, draw, scale=.6] (0) at (0,0){ab};
        \node[circle, draw, scale=.6] (1) at (0.5,1){cd};
        \node[circle, draw, scale=.6] (2) at (-0.5,1){ef};
        \draw (0)--(1);
        \draw (1)--(2);
        \draw (0)--(2);

        \node[circle, draw, scale=.6] (0) at (3,0){ab};
        \node[circle, draw, scale=.6] (1) at (3.5,1){cd};
        \node[circle, draw, scale=.6] (2) at (2.5,1){ed};
        \draw (0)--(1);
        \draw (1)--(2);
        \draw (0)--(2);

        \node[circle, draw, scale=.6] (0) at (6,0){ab};
        \node[circle, draw, scale=.6] (1) at (6.5,1){cd};
        \node[circle, draw, scale=.6] (2) at (5.5,1){ad};
        \draw (0)--(1);
        \draw (1)--(2);
        \draw (0)--(2);

        \node[circle, draw, scale=.6] (0) at (9,0){ab};
        \node[circle, draw, scale=.6] (1) at (9.5,1){bd};
        \node[circle, draw, scale=.6] (2) at (8.5,1){ad};
        \draw (0)--(1);
        \draw (1)--(2);
        \draw (0)--(2);
    \end{tikzpicture}
    \end{center}

In the first case we obtain that $n=6$, in the second case we obtain that the vertex $d$ has degree at least three since $(d,c), (d,c) \in E(G)$ and $ a$ or $b$ is a also a neighbour of $d$ which is a contradiction. In the third case we obtain $n=4$ since $(c,d),(a,d),(a,b),(b,c)\in E(G)$, and in the fourth case we have $n=3.$ Hence, $F_2^2(C_n)$ is triangle free for $n \neq 3,4,6$ thus $\omega (F_2^2(C_n))=2.$

Now consider  the chromatic number. If $n$ is odd, by Theorem \ref{prod} $F_2^2(C_n)=C_n \square P_k$ thus $\chi (F_2^2(C_n)) = max \{ \chi(C_n), \chi( P_k)\}=3$.

 If $n=2k$, by Proposition \ref{comp}, $F_2^2(C_n)$ contains one bipartite connected component, thus we only need to show that the non-bipartite connected component can be coloured with three colours. 
 
Recall the description of the non-bipartite component of $F_2^2(C_n)$ given in Remark \ref{nonbi}. We are going to colour $V(F_2^2(C_n))$ as follows. Let $c: V(F_2^2(C_n)) \rightarrow \{0,1,2\}$ given by $c(x,y)= i $ if $x<y$ and $x\equiv i (\text{ mod }3)$ with the convention that $n$ is going to be considered smaller than $1.$ It is easy to see that this is a proper colouring if $n \not\equiv 1 (\text{mod }3).$ When $n \equiv 1(\text{mod }3)$, we consider the colouring $c':V(F_2^2(C_n)) \rightarrow \{0,1,2\}$ given by
\begin{align*}
   & c'(x,y)= i  \text{ if } x<y, x>3 \text{ and } x\equiv i (\text{ mod }3) \\
   &c'(1,x)=2 \text{ for every } 2 \leq x < n  \text{ and } c'(2,x)=1 \text{ for } 2<x\leq n.
\end{align*}
  Then $c'$ is a proper colouring of the non bipartite component of $F_2^2(C_{2k})$ with three colours.
\end{proof}

\begin{teo}
   For any even integer $n,$ the independence number of $F_2^2(C_n)$ is $$\alpha(F_2^2(C_n))=
       \frac{n}{8}(n-2) + \binom{\frac{n}{2}}{2} .$$
\end{teo}

\begin{proof}
    By Proposition \ref{comp}, $F_2^2(C_n)$ has two connected components $A$ and $B$, such that $B$ is bipartite. It is easy to see that $\alpha(B)=\binom{\frac{n}{2}}{2}.$ Consider the sets $V_i$ described in Remark \ref{nonbi} and notice that $B\setminus V_{k'} \cong C_n\square P_{\frac{k'-1}{2}}$ and $\alpha(B)\leq \alpha (B\setminus V_{k'})+ \alpha(V_{k'})= \frac{n}{2}(\frac{k'-1}{2})+\alpha(V_{k'})$. If $n \equiv 0\text{ (mod $4$)}$ then $\alpha(V_{k'})=\frac{n}{2}-1=k'$ and if $n \equiv 2\text{ (mod $4$)}$ then  $\alpha(V_{k'})=\lfloor \frac{n}{4} \rfloor$ and $k'= \frac{n}{2}.$
    
    However, due to the edges between $V_{k'}$ and $V_{k'-2}$, we can see that if $n \equiv 0\text{ (mod $4$)}$ then $\alpha(B)=(\frac{k'-1}{2})(\frac{n}{2})+\frac{n}{4}=\frac{n}{4}k'=\frac{n}{8}(n-2)$ and if $n \equiv 2\text{ (mod $4$)}$ then $\alpha(B)=\alpha(B\setminus V_{k'})= (\frac{k'-1}{2})(\frac{n}{2})= (\frac{\frac{n}2-1}{2})\frac{n}{2}=\frac{n}{8}(n-2).$
\end{proof}

\begin{coro}
Let $n$ be an even integer. Then the domination and independent domination numbers of $F_2^2(C_n)$ are bounded by $\gamma (F_2^2(C_n))\leq i(F_2^2(C_n)) \leq   \frac{n}{8}(n-2) + \binom{\frac{n}{2}}{2}.$
\end{coro}

\begin{figure}[h!]
    \centering
    \includegraphics[scale=0.4]{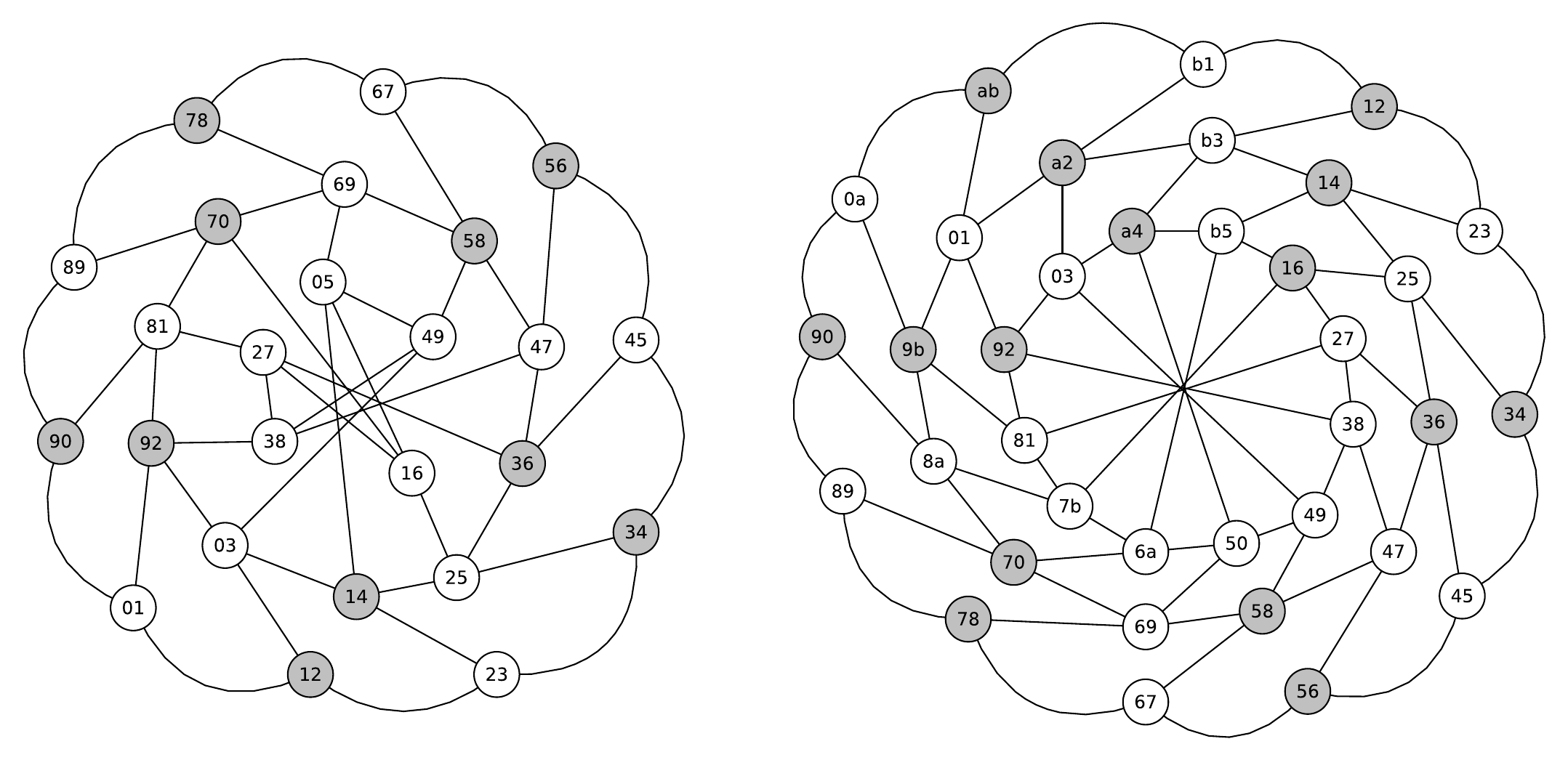}
    \caption{The non bipartite components of the graphs $F_2^2(C_{10})$ and $F_2^2(C_{12})$ respectively, where the gray vertices represent a maximal independent set.}
    \label{flores}
\end{figure}

In the following theorem we give  conditions on the graph $G$ for the generalized token graph $F_2^2(G)$ to be connected.  
 
\begin{figure}[h!]
\begin{center}
		\begin{tikzpicture}
		\tikzstyle{every node}=[minimum width=0pt, inner sep=2pt, circle]
			\draw (-2,0) node[draw] (0) {};
			\draw (-0.6180339887498951,-1.9021130325903073) node[draw] (1) {};
			\draw (1.618033988749894,-1.1755705045849465) node[draw] (2) {};
			\draw (1.6180339887498951,1.1755705045849458) node[draw] (3) {};
			\draw (-0.6180339887498946,1.9021130325903073) node[draw] (4) {};
			\draw (-2.7,0) node[draw] (5) { $K_{2,m_k}$};
			\draw (-0.69,-2.6) node[draw] (6) {$K_{2,m_4}$};
			\draw (2.3,-1.31) node[draw] (7) {$K_{2,m_3}$};
			\draw (2.3,1.35) node[draw] (8) {$K_{2,m_2}$};
			\draw (-0.69,2.6) node[draw] (9) {$K_{2,m_1}$};
			\draw[dotted]  (0) edge (1);
			\draw  (1) edge (2);
			\draw  (2) edge (3);
			\draw  (3) edge (4);
			\draw   (0) edge (4);
		\end{tikzpicture}
	\end{center} \caption{A graph obtained from an odd cycle by glueing complete bipartite graphs to the vertices of the cycle.}
    \label{ciclo5kmn}
\end{figure}
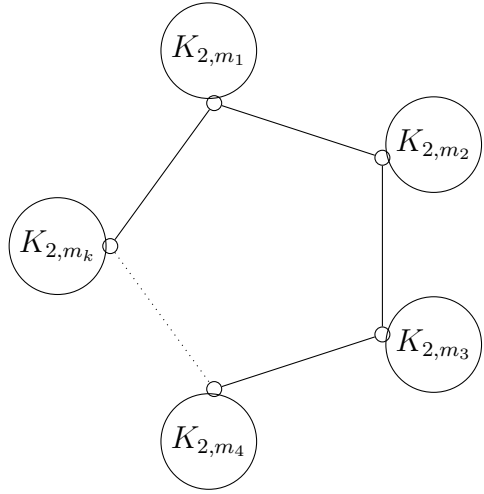

\begin{teo}\label{conex}
    Let $G$ be a non-bipartite graph. Then,
    \begin{itemize}
        \item If $G$ does not contain a set of two or more leaves which have a common neighbour, then $F_2^2(G)$ is connected, and
        \item if $G$ contains a set of $k\geq2$ leaves which have a common neighbour, then $F_2^2(G)$ consists of a single non-trivial connected component together with $\binom{k}{2}$ isolated vertices.
    \end{itemize}
\end{teo}
\begin{proof}
    Proceed by induction over $|V(G)|.$ If $|V(G)|=3$ then $G$ is either a cycle of length three and thus $F_2^2(G)$ is a cycle of length three, or $G$ is a path with three vertices in which case $F_2^2(G)$ is the disjoint union $K_2\cup K_1.$ 
    Now assume that $|V(G)|=n$. Consider first the case when $G$ does contain a set of two or more leaves which have a common neighbour and let $v$ be a vertex such that $G\setminus \{v\}$ is not bipartite (if $G\setminus{v}$ is bipartite for every $v\in V(G)$ then $G$ must be an odd cycle thus by \ref{prod}, $F_2^2(G)$ is connected) 
    and does not contain a set of two or more leaves which have a common neighbour. Such a vertex exists unless $G$ is either a cycle, or, a graph obtained from a single odd cycle by gluing complete bipartite graphs $K_{2,m}$ to vertices of the odd cycle as in Figure \ref{ciclo5kmn}.
    
    Assume first such vertex $v$ exists. Then $F_2^2(G\setminus v)$ is connected, and we only need to show that there exists a path from a vertex of the form $vw$ to a vertex in $F_2^2(G\setminus v)$. Let $u$ be a neighbour of $v$ and $z$ be a neighbour of $w$. If $z \neq v$ then $(vw, uz)\in E(F_2^2(G))$. If $z=v$ then $(vw, vu)\in E(F_2^2(G))$. Then either $(w,u)\in E(G)$ in which case both $vu$ and $vw$ are adjacent to $wu$, or there exists a neighbour $s$ of $u$ different from $w$ and $v$ in which case $(vu, su)\in E(F_2^2(G))$ and hence $(F_2^2(G))$ is connected.
    
Now, assume that $G\setminus\{v\}$ is not bipartite but contains an independent set of vertices $\{x_1, \dots, x_m\}$ which have as common neighbour a vertex $w.$ Then by induction hypothesis, $F_2^2(G)$ has a single non-trivial connected component and $\binom{m}{2}$ isolated vertices of the form $x_ix_j.$ Then, in $F_2^2(G)$, the vertices $x_ix_j$ are adjacent to the vertex $vw,$ thus they are no longer isolated. As in the previous case, $vw$ must be connected to the non trivial connected component of $F_2^2(G)$ and hence, $F_2^2(G)$ is connected.

Assume that $G$ is not bipartite and contains an independent set of vertices $\{x_1, \dots, x_m\}$ which have as common neighbour a vertex $w.$ We remove the vertex $x_m$, thus, if $m>2$, by induction hypothesis, $F_2^2(G\setminus x_m)$ has a single non trivial connected component and isolated vertices. Then in $F_2^2(G)$, the vertex $x_mx_i$ is isolated for every $ 1 \leq i \leq m-1$. Consider $x_mv$ for some $v\neq x_i$, and let $u$ be a neighbour of $v$ with $u\neq w.$ Then $(x_mv,wu)\in E(F_2^2(G)).$
Finally, if $m=2$, by induction hypothesis $F_2^2(G\setminus x_m)$ is connected and $x_1x_2$ is an isolated vertex in $F_2^2(G).$ As in the previous case, any vertex of the form $x_2v$ must be adjacent to a vertex of $F_2^2(G\setminus x_m).$
\end{proof}

\begin{lema}\label{prod}
    Given a disjoint union of graphs $G \sqcup H$, $F_2^2(G\sqcup H)= F_2^2(G) \sqcup F_2^2(H) \sqcup G \times H.$
\end{lema}

\begin{proof}
    Consider the sets $X,Y,Z\subset V(F_2^2(G \sqcup H))$ given by $X=\{(x,g): x,g \in V(G)\}, Y=\{(y,h): y,h \in V(H)\}$ and $Z=\{(g,h): g\in V(G), h\in V(H)\}.$ Then $ V(F_2^2(G \sqcup H))=X\sqcup Y\sqcup Z$ and notice that since there are no edges between $G$ and $H$, there are no edges between $X, Y$ and $Z.$ The subgraphs induced by $X$ and $Y$ are precisely $F_2^2(G)$ and $F_2^2(H)$ respectively. The subgraph spanned by $Z$ has $V(\langle Z \rangle)=Z=V(G)\times V(H)$ and $E(\langle Z \rangle)=\{((g,h), (z,w)): (g,z)\in E(G), (h,w)\in E(H)\}$ thus $\langle Z \rangle = G \times H$ and hence $F_2^2(G\sqcup H)= F_2^2(G) \sqcup F_2^2(H) \sqcup G \times H.$
\end{proof}

We have the following characterization for graphs $G$ such that its generalized token graph $F_2^2(G)$ is a bipartite graph.
\begin{teo}
    The graph $F_2^2(G)$ is bipartite if and only if $G=\bigsqcup_{i=1}^nP_{m_i}$.
\end{teo}
\begin{proof}
    Assume that $F_2^2(G)$ is bipartite and notice that if $G$ contains a vertex of degree three or higher, then $F_2^2(G)$ contains a subgraph isomorphic to $K_3.$ This means that $G$ is either a path or a cycle. If $G$ is an odd cycle, then $F_2^2(G)$ contains an odd cycle. Thus assume that $G=C_{2n}$ with vertex set $\{1, 2, \dots, 2n\}$ and edges $\{(i,i+1), (2n,1)\}$ for $1\leq i \leq 2n-1$ with $n\geq2.$   If $n$ is even, then the vertices $\{1,n\},\{2,n+1\}, \dots, \{n+1,2n\}$ form a cycle of length $n+1$ which is odd. If $n$ is odd, then the vertices $\{1,n+1\},\{2,n+2\},\dots, \{n,2n\}$ form a cycle of length $n.$
    Hence $G$ cannot contain a cycle and must be a union of paths.
    Now assume that $G$ is a union of paths, $G =\bigcup_{i=1}^m P_{m_i}$. Then since, $F_2^2(P_{m_i})$ and $P_{m_i}\times P_{m_j}$ are bipartite, by Lemma \ref{prod}, the graph $F_2^2(G)$ is bipartite.
\end{proof} 

\begin{coro}
   The chromatic number $\chi(F_2^2(G))=2$ if and only if $G=\bigsqcup_{i=1}^nP_{m_i}$. 
\end{coro}
 \begin{coro}
      Then the clique number  $\omega(F_2^2(G))=2$ if and only if $G=\bigsqcup_{i=1}^nP_{m_i}$.
 \end{coro}

\begin{teo}
    The independence number of $F_2^2(P_n)$ is $$\alpha(F_2^2(P_n))=\left\lfloor \frac{n}{2} \right\rfloor \left\lceil \frac{n}{2} \right\rceil.$$
\end{teo}
\begin{proof}
    Assume that $V(P_n)=\{1, \dots, n\}$ and $E(P_n)=\{(i,i+1)\}_{i=1}^n.$ Consider the set of vertices $A=\{(1,2),\dots, (1,n), (3,4), \dots, (3,n), \dots\}\subset V(F_2^2(P_n))$ where  the last vertex of $A$ is $(n,n-1)$ if $n$ is even and, if $n$ is odd then the last two vertices of $A$ are $(n-1,n-2)$ and $(n, n-2).$
    It is easy to see that $A$ is an independent set, thus $\alpha(F_2^2(P_n))\geq |A|.$ Notice that $|A|=n-1 + n-3 + \dots + f(n)$ where $f(n)= 1$ if $n $ is even and $f(n)=2$ if $n$ is odd. This means that $$\alpha(F_2^2(P_n)) \leq |A|=\left\{\begin{array}{cl} k^2 & \text{ if }n=2k \\
       k(k+1) & \text{ if }n=2k+1
    \end{array}\right. .$$
    
    On the other hand, notice that $F_2^2(P_n)$ contains $n-1$ vertex disjoint paths of lengths $n-1, n-2, \dots, 1$ thus $\alpha(F_2^2(P_n))\leq \lfloor \frac{n-1}{2} \rfloor + \lfloor \frac{n-2}{2} \rfloor + \dots 1=x$ 
    
    If $n$ is odd then $\alpha(F_2^2(P_n)) \leq x = 2(\sum\limits_{i=1}^{\frac{n-1}{2}} i)= \frac{n-1}{2}\frac{n+1}{2}.$ If $n$ is even, then $\alpha(F_2^2(P_n)) \leq x = \frac{n}{2} + 2(\sum\limits_{i=1}^{\frac{n-2}{2}} i)= (\frac{n}{2})^2. $
    Hence, $x=|A|$ and $\alpha(F_2^2(P_n))=\left\lfloor \frac{n}{2} \right\rfloor \left\lceil \frac{n}{2} \right\rceil.$
\end{proof}

In \cite{I}, the following values for the independence numbers of direct products of graphs were obtained.
\begin{teo}\cite{I}
    The independence number of the direct product of paths is $\alpha (P_m\times P_n)=\frac{mn}{2}$ if both $m$ and $n$ are even, and $\alpha (P_m\times P_n)=\frac{m(n+1)}{2}$  if $m$ is odd.
\end{teo}

\begin{coro}
    Let $P_m$ and $P_n$ be two paths. Then $$\alpha(F_2^2(P_n\sqcup P_m))=\left\lfloor \frac{n}{2} \right\rfloor \left\lceil \frac{n}{2} \right\rceil +\left\lfloor \frac{m}{2} \right\rfloor \left\lceil \frac{m}{2} \right\rceil + \alpha (P_n \times P_m) $$ where $\alpha (P_m\times P_n)=\frac{mn}{2}$ if both $m$ and $n$ are even, and $\alpha (P_m\times P_n)=\frac{m(n+1)}{2}$  if $m$ is odd.
\end{coro}

We conclude this section with the following result.

\begin{teo}\label{comp}
Let $G$ be a graph. Then $\bigcup_{i=1}^n A(F_n^i(G))$ is the set of arcs of a complete graph if and only if  the complement $G^c$ of $G$ does not contain as a subgraph any complete graph $K_{m+k}$ or complete bipartite graph $K_{m,k}$ for some $m,k$ with $m+k>n$.
\end{teo}
\begin{proof}
    Assume that $\bigcup_{i=1}^n A(F_n^i(G))$ is the set of arcs of a complete graph and suppose by contradiction that $G^c$ contains a subgraph isomorphic to $K_{m,k}$ with $m+k>n.$ Assume the vertices of the $K_{m,k}$ are $\{x_1, \dots, x_k\}$ and $\{y_1, \dots, y_m \}$ and let $A=\{x_1, \dots, x_k, x_{k+1},\dots, x_n\}$ and $B=\{y_1, \dots, y_m, y_{m+1},\dots, y_n\}$. Since $(A,B)$ is an edge in $\bigcup_{i=1}^n A(F_n^i(G))$, each vertex of $\{x_1, \dots, x_k\}$ must be adjacent to and paired up with a vertex of $\{y_{m+1},\dots, y_n\}$ thus in particular $k \leq n-m$ but this implies that $m+k \leq n $ which is a contradiction.

Now assume again by contradiction, that $\bigcup_{i=1}^n A(F_n^i(G))$ is not the set of arcs of a complete graph. This is, there are vertices $A=\{x_1,\dots, x_n\}$ and $ B=\{y_1, \dots, y_n\}$ such that $(A,B)\notin \bigcup_{i=1}^n A(F_n^i(G))$. This means that we cannot pair up the vertices of $A$ with the vertices of $B$ as $\{x_i,y_i\}$ such that either $x_i=y_i$ or $(x_i, y_i)\in E(G).$ In particular this means that there exists a subset $X \subset A$ of $k$ vertices with $1\leq k \leq n $ such that the intersection of its neighbourhood with $B$, $N(X)\cap B$ consists of less than $k$ vertices. Then, since there are no vertices from $X$ to $B\setminus (N(X)\cap B)$ we have that $X$ together with $B\setminus (N(X)\cap B)$ form a subgraph in $G^c$ which contains a complete bipartite subgraph $K_{k,m}$ with $m\geq n-k+1$ a contradiction.
\end{proof}

\begin{coro}
    If $G^c$ does not contain as a subgraph $K_3$ or $K_{1,2}$, then $F_2^2(G)= (F_2(G))^c.$
\end{coro}

\subsection{On the automorphism group of $F_n^k(G).$}
In this last section of the paper, we give a brief analysis of the autormorphism group of the generalized token graph.

\begin{teo}
    Given a graph $G,$ $Aut(G)\leq Aut(F_n^k(G)).$
\end{teo}
\begin{proof}
    Let $f\in Aut(G)$ and define  $\phi_f:F_n^k(G)\rightarrow F_n^k(G)$ given by $\phi(x_1,\dots, x_n)= (f(x_1),\dots, f(x_n)).$ We claim that $\phi_f \in Aut(F_n^k(G)).$
    Let $((x_1,\dots,x_n), (y_1,\dots, y_n))\in E(F_n^k(G))$  and assume without loss of generality that $(x_i,y_i)\in E(G)$ for $1\leq i \leq k$ and $x_j=y_j$ for $j>k.$ Then, since $f\in Aut(G)$, we have $(f(x_i), f(y_i))\in E(G)$  for $1 \leq i \leq k$, thus $((f(x_1), \dots, f(x_n),(f(y_1),\dots,f(y_n))\in E(F_n^k(G)).$
\end{proof}

\begin{teo}\label{autmn}
    For $m,n>1$ we have $Aut(F_2^2(K_{m,n}))=\mathcal{S}_{mn}\times \mathcal{S}_{\binom{n}{2}}\times \mathcal{S}_{\binom{m}{2}}.$
\end{teo}
\begin{proof}
    By Theorem \ref{kmn},  $F_2^2(K_{m,n})$ is isomorphic to the disjoint union $K_{mn}\sqcup K_{\binom{n}{2},\binom{m}{2}}$ thus $$Aut(F_2^2(K_{m,n}))=Aut(K_{mn})\times Aut(K_{\binom{n}{2},\binom{m}{2}})= \mathcal{S}_{mn}\times (\mathcal{S}_{\binom{n}{2}}\times \mathcal{S}_{\binom{m}{2}}).$$
\end{proof}

Notice that Theorem \ref{autmn} does not hold when either $m$ or $n$ is equal to one since then the generalized token graph contains several isolated vertices which can be permuted with each other by an automorphism.

In \cite{hik}, the following theorem was proven.
\begin{teo}[Theorem 6.13,\cite{hik}]
    The automorphism group of the Cartesian product of connected prime
graphs is isomorphic to the automorphism group of the disjoint union of the factors.
\end{teo}

\begin{teo} \label{autprod}
    For any odd integer $n>3$, $Aut(F_2^2(C_n))=D_n\times \mathcal{S}_2.$
\end{teo}
\begin{proof}
    By Theorem \ref{prod}, we know that $F_2^2(C_n)$ is isomorphic to the product $C_n\square P_{\frac{n-1}{2}}.$ By Theorem \ref{autprod}, $Aut(F_2^2(C_n))=Aut(C_n)\times Aut(P_{\frac{n-1}{2}})=D_n \times \mathcal{S}_2.$
\end{proof}

\begin{ejem}
    Consider the diamond graph $D=K_4\setminus e$ for $e\in E(K_4).$ Then $Aut(D)=\mathcal{S}_2\times \mathcal{S}_2$ and $F_2^2(D)$ is the graph on six vertices which has one leaf, one vertex of degree five, and every other vertex has degree four thus $Aut(F_2^2(D))=\mathcal{S}_4$ as we can see in Figure \ref{diamond}. 
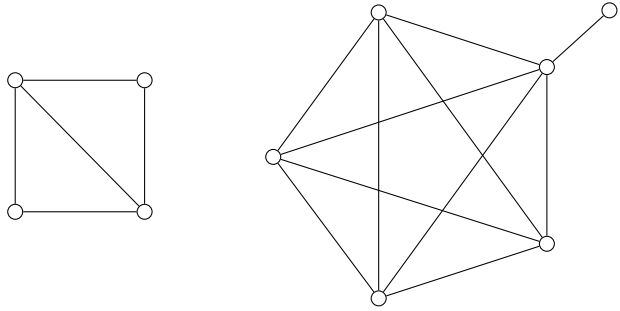
\begin{figure}
    \begin{center}
		\begin{tikzpicture}
		\tikzstyle{every node}=[minimum width=0pt, inner sep=2pt, circle]
			\draw (-2.02,-0.02) node[draw] (0) {};
			\draw (-0.6180339887498951,-1.9021130325903073) node[draw] (1) { };
			\draw (1.618033988749894,-1.1755705045849465) node[draw] (2) { };
			\draw (1.6180339887498951,1.1755705045849458) node[draw] (3) {} ;
			\draw (-0.6180339887498946,1.9021130325903073) node[draw] (4) { };
			\draw (-3.73,1) node[draw] (5) { };
			\draw (-5.45,-0.75) node[draw] (6) {};
			\draw (-3.73,-0.75) node[draw] (7) {};
			\draw (-5.45,1) node[draw] (8) { };
			\draw (2.45,1.93) node[draw] (9) { };
			\draw  (0) edge (1);
			\draw  (0) edge (2);
			\draw  (0) edge (3);
			\draw  (0) edge (4);
			\draw  (1) edge (2);
			\draw  (1) edge (3);
			\draw  (1) edge (4);
			\draw  (2) edge (3);
			\draw  (2) edge (4);
			\draw  (3) edge (4);
			\draw  (5) edge (7);
			\draw  (5) edge (8);
			\draw  (6) edge (7);
			\draw  (6) edge (8);
			\draw  (7) edge (8);
			\draw  (3) edge (9);
		\end{tikzpicture}
	\end{center}
    \caption{The diamond graph $D$ and $F_2^2(D).$} \label{diamond}
    \end{figure}
\end{ejem}

It is well known that, for usual token graphs, if a graph $G$ does not contain induced cycles of length four, or induced diamonds ($K_4\setminus e)$ then $Aut(G)=Aut(F_k(G))$ for every $k \neq \lfloor \frac{n}{2} \rfloor.$ 
So, we conclude with the following open question.

\begin{conje}
    For which graphs $G$ does the following equality hold?
    $$Aut(G)=Aut(F_2^2)$$
\end{conje}


\section*{Statements and Declarations}
The authors have no relevant financial or non-financial interests to disclose.
Data sharing is not applicable to this article, as no new data was created or analyzed in this study.

\end{document}